\documentclass[10pt]{article}
\usepackage[utf8]{inputenc}
\usepackage{amsfonts,amsmath,amssymb,amsthm}
\usepackage{graphicx,verbatim}
\usepackage{float}
\usepackage{subfig}
\usepackage{xcolor}
\usepackage{hyperref}

\RequirePackage{geometry}
\newcommand{\pgap}{\vspace{0.3cm}}

\newcommand{\bR}{\mathbb{R}}
\newcommand{\bE}{\mathbb{E}}

\newcommand{\bI}{\mathbf{1}}
\newcommand{\bP}{\mathbb{P}}

\newcommand{\bN}{\mathbb{N}}
\newcommand{\bC}{\mathbb{C}}

\newcommand{\cov}{\operatorname{Cov}}

\newcommand{\td}{\text{d}}

\newcommand{\FGF}{\mathrm{FGF}}
\renewcommand{\phi}{\varphi}

\renewcommand{\epsilon}{\varepsilon}

\newtheorem{proposition}{Proposition}[section]
\newtheorem{theorem}[proposition]{Theorem}

\newtheorem{remark}[proposition]{Remark}

\title{Excursion Fluctuations and Spectral Universality in Gaussian Fields}
\author{Dmitry Beliaev \footnote{Mathematical Institute, University of Oxford, UK } \and Akshay Hegde \footnote{Department of Mathematics, National University of Singapore, Singapore}}
\date{}

\begin{document}

\maketitle

\begin{abstract}
We study the large-scale spatial fluctuations of excursion volumes for a class of smooth stationary Gaussian fields. In the case of Berry's random wave model in dimension $d \geq 2$, we show that the spatial fluctuations for fixed $u>0$ converge to the fractional Gaussian field $(-\Delta)^{-1/4}W$ in the space of tempered distributions $\mathcal S'(\bR^d)$, where $W$ is the $d$-dimensional Gaussian white noise. This explains the long-range correlations in the apparent filament structure of the Random Plane Wave model. For a class of smooth planar Gaussian fields whose spectral density has a power-law singularity at the origin, we prove convergence to fractional Gaussian fields with an index determined by the singularity exponent. More generally, the results illustrate that, for stationary random measures, large-scale spatial fluctuations are determined by the behaviour of the spectral measure density exponent near zero. 
\end{abstract}

\setcounter{tocdepth}{1}


\section{Introduction}

\begin{figure}[ht]
    \centering
    \includegraphics[width=0.9\linewidth]{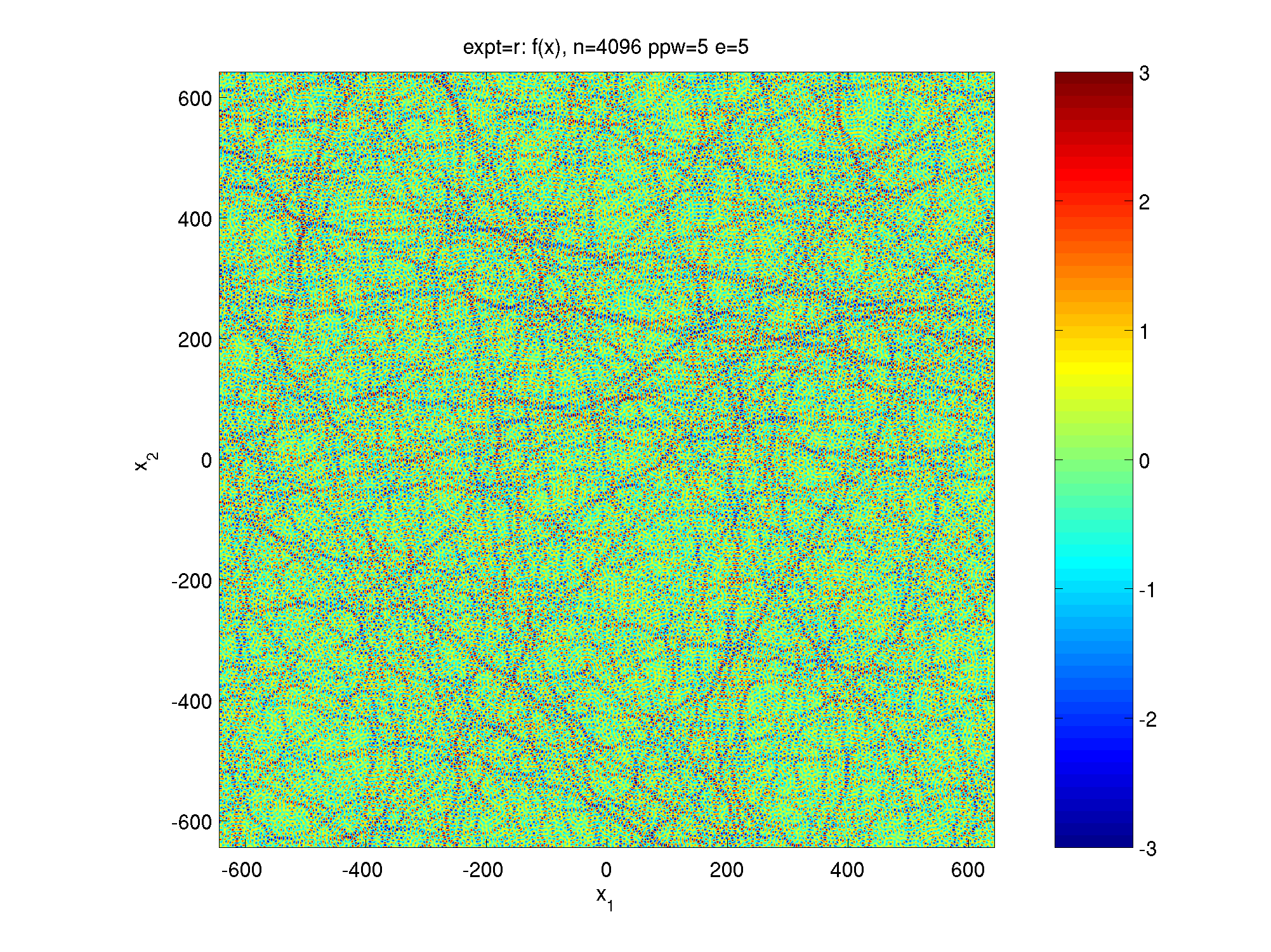}
    \caption{Heatmap of a sample of the Random Plane Wave model. Picture credit: Alex Barnett \cite{Warms_movie}}
    \label{fig:RPW-filament}
\end{figure}

The original motivation for this work was to understand the filaments: apparent line-like patterns visible in samples of the random plane wave (See Figure \ref{fig:RPW-filament} for an example that initially sparked our curiosity.) Our earlier work \cite{beliaev_hegde} on high local maxima showed that, under very general assumptions, the point process corresponding to high-level local maxima has a Poisson scaling limit, so filamentary structure cannot be detected at the level of sparse extreme points. The present paper does not give a mathematical definition or direct proof of filament geometry. Instead, it identifies a different mechanism: long-range spatial correlations can appear at a \emph{fixed} non-zero level, and their scaling limit is governed by the behaviour of the spectral measure near the origin. In this sense, our results help reconcile the apparent visual structure with the Poissonian behaviour of high maxima by showing that the two phenomena belong to different observables and different scaling regimes.

\pgap 

\pgap 

Our novel observation is the following. Let $M$ be a centered stationary random measure on $\bR^d$ and assume its
covariance measure has spectral\footnote{Throughout this paper we use the Fourier transform convention $
    \widehat{g}(\xi)=\int_{\bR^2} e^{-i x \cdot \xi}g(x) \td x .
$} 
density $m$:
\[
    \mathrm{Cov}(M(\varphi),M(\psi))
    =\frac{1}{(2\pi)^d}\int_{\bR^d}
    \widehat\varphi(\xi)\widehat\psi(-\xi)m(\xi)\td \xi .
\]
For $\varphi_R(x):=\varphi(x/R)$,
\[
    \mathrm{Cov}(M(\varphi_R),M(\psi_R))
    =
    \frac{R^d}{(2\pi)^d}\int_{\bR^d}
    \widehat\varphi(\eta)\widehat\psi(-\eta)m(\eta/R)\td \eta .
\]
Assume
\[
    m(\xi)\sim c|\xi|^{-2s}\quad (\xi\to0),
    \qquad c>0,\qquad s<d/2.
\]
Then we should have\footnote{Here we ignore exceptional cases where the Hurst parameter $H=s-d/2$ is an integrer or outside of a certain range. The relevant ranges will be clarified later on. } 
\[
\begin{aligned}
    R^{-d-2s}\mathrm{Cov}(M(\varphi_R),M(\psi_R))
    \to &
    \frac{c}{(2\pi)^d}\int_{\bR^d}
    \widehat\varphi(\eta)\widehat\psi(-\eta)|\eta|^{-2s}\td \eta 
    \\ 
    =&c \int_{\bR^2}\int_{\bR^2} \varphi(x)\psi(y)|x-y|^{2s-d} \td x \td y .
\end{aligned}    
\]

Hence the natural normalisation is
\[
    X_R(\varphi):=R^{-d/2-s}\{M(\varphi_R)-\bE[M(\varphi_R)]\}.
\]
If the characteristic functionals of $X_R$ converge to the Gaussian
functional with the above covariance, and if the laws of $X_R$ are tight in
$\mathcal{S}'(\bR^d)$, then
\[
    X_R\Rightarrow \sqrt c\,\mathrm{FGF}_s(\bR^d)
    :=\sqrt c\,(-\Delta)^{-s/2}W ,
\]
where $W$ is the $L^2(\bR^d)$ white noise. 
Here $\mathrm{FGF}_s$ is the fractional Gaussian field. The right-hand side gives an informal definition. FGF is not a function but a random distribution which can be defined as a Gaussian field indexed by test functions with covariance given by 
\[          
\mathrm{Cov}\big(\mathrm{FGF}_s(\varphi),\mathrm{FGF}_s(\psi)\big)
    =
    \frac{1}{(2\pi)^d}\int_{\bR^d}
    \widehat\varphi(\xi)\widehat\psi(-\xi)|\xi|^{-2s}\td \xi 
\]
where $\phi$ and $\psi$ are suitable test functions. This is exactly the covariance given by the heuristics above. For a comprehensive survey of the fractional Gaussian fields, we refer the readers to \cite{lodhia_fractional_2016}. A short introduction will be given in Appendix \ref{ss: FGF}.

When this paper was essentially completed, the preprint \cite{gass2026} was posted on the ArXiv. It studies closely related fluctuation phenomena for Berry’s random wave model. Our work was carried out independently. There is some overlap in the random-wave case: in particular, our Theorem 2.1 is essentially a partial case of \cite[Proposition 2.4.3]{gass2026}. The emphasis of the two papers is, however, different. The focus of \cite{gass2026} is on the generality of observables in the random-wave setting, whereas here the focus is on the generality of fields and on the dependence of the universality class on the low-frequency singularity of the spectrum. In particular, beyond the $\FGF_{1/2}$ limit for Berry’s model, we obtain a family of different fractional-Gaussian limits for planar fields with spectral density behaving like $|\xi|^{-\alpha}$ near the origin.

\subsubsection*{Acknowledgement }
The authors would like to thank the Isaac Newton Institute for Mathematical Sciences, Cambridge, for support and hospitality during the programme Geometric spectral theory and applications, where part of this work on this paper was undertaken. This work was supported by EPSRC grant EP/Z000580/1.

\section{Main results}

\subsection{Berry's Random Wave Model}

 For $d\geq 2$, let
\[
    \omega_{d-1}:=\operatorname{area}(S^{d-1})
    =\frac{2\pi^{d/2}}{\Gamma(d/2)},
    \qquad
    \mu_d:=\frac{1}{\omega_{d-1}}\sigma_{S^{d-1}},
\]
where $S^{d-1} \subset \bR^d$ is the unit sphere, $\sigma_{S^{d-1}}$ is the uniform hypersurface measure on $S^{d-1}$, and $\Gamma$ is the Gamma function. Let $F_d:\bR^d\to\bR$ be the stationary centered Gaussian field with
spectral measure $\mu_d$. Thus
\[
    K_d(x-y):=\bE[F_d(x)F_d(y)]
    =\int_{S^{d-1}}e^{i(x-y)\cdot\theta}\td\mu_d(\theta). \]
In dimension $2$ this is $K_2(z)=J_0(|z|)$ ($J_0$ is the zeroth Bessel function of first kind ), while in dimension $3$ it is $K_3(z)=\sin |z|/|z|$.

\pgap 

The main object of study in this section is the Berry's random wave model $F_d$ and in particular, the random plane wave (RPW) field, the two-dimensional case. Figure \ref{fig:RPW-filament} shows the heatmap of a sample of RPW. 

\pgap 

Our goal is to understand the scaling behaviour of excursion sets $\{x: F_d(x)> u\}$.  In the terminology of the Introduction, we are studying the random measure $\bI[F_d(x)>u] \td x$ and its linear functionals. We are going to work with the following setting. For fixed $u>0$, $R>0$, and $\varphi\in\mathcal S(\bR^d)$ define the linear statistics of excursion volume
\begin{equation} \label{eq:area-statistics}
    A^{(d)}_{R,u}(\varphi)
    :=
    \int_{\bR^d}\varphi(x/R)\bI[F_d(x)>u]\td x.
\end{equation}

Let $h_d=(-\Delta)^{-1/4}W$ be the fractional Gaussian field  $\FGF_{1/2}$ on $\bR^d$ with
covariance
\begin{equation} \label{eq:monochromatic-fgf-covariance}
    \bE[h_d(\varphi)h_d(\psi)]
    =
    \frac{1}{(2\pi)^d}
    \int_{\bR^d}
    \frac{\widehat\varphi(\xi)\widehat\psi(-\xi)}{|\xi|}
    \td \xi .
\end{equation}

\begin{theorem}
\label{thm:main-thm}
Fix $d\geq2$ and $u>0$. As $R\to\infty$,
\[
    R^{-(d+1)/2}
    \big(A^{(d)}_{R,u}-\bE[A^{(d)}_{R,u}]\big)
    \Rightarrow
    c_{d,u}h_d
    \qquad \text{in } \mathcal S'(\bR^d),
\]
where
\[
    c_{d,u}
    =
    u\gamma(u)
    \left(
        2^{d-2}\pi^{(d-1)/2}
        \frac{\Gamma(d/2)^2}{\Gamma((d-1)/2)}
    \right)^{1/2},
    \qquad
    \gamma(u)=\frac{1}{\sqrt{2\pi}}e^{-u^2/2}.
\]
\end{theorem}

\begin{remark}
Note that the parameter $s=1/2$ is independent of the dimension. Somewhat more conventional Hurst exponent $H=s-d/2$ does depend on the dimension.     
\end{remark}

This theorem is the case of $X_4$ in \cite[Proposition 2.4.3]{gass2026}. There is a small difference in formulation. We do not have exceptional values of $u$, but when $u=0$, the constant $c_u$ vanishes, and the theorem no longer implies convergence to FGF. This means that in this case $R^{(d+1)/2}$ is not the right scaling. This is yet another manifestation of Berry's cancellation at the level $u=0$. 

\pgap 

Another important point is that here, unlike previous studies aimed at understanding the filaments, we do not send $u$ to infinity. We claim that the long-range correlations appear at any \emph{fixed}  non-zero level. 

\begin{figure}[ht]
    \centering
    \includegraphics[width=0.49\linewidth]{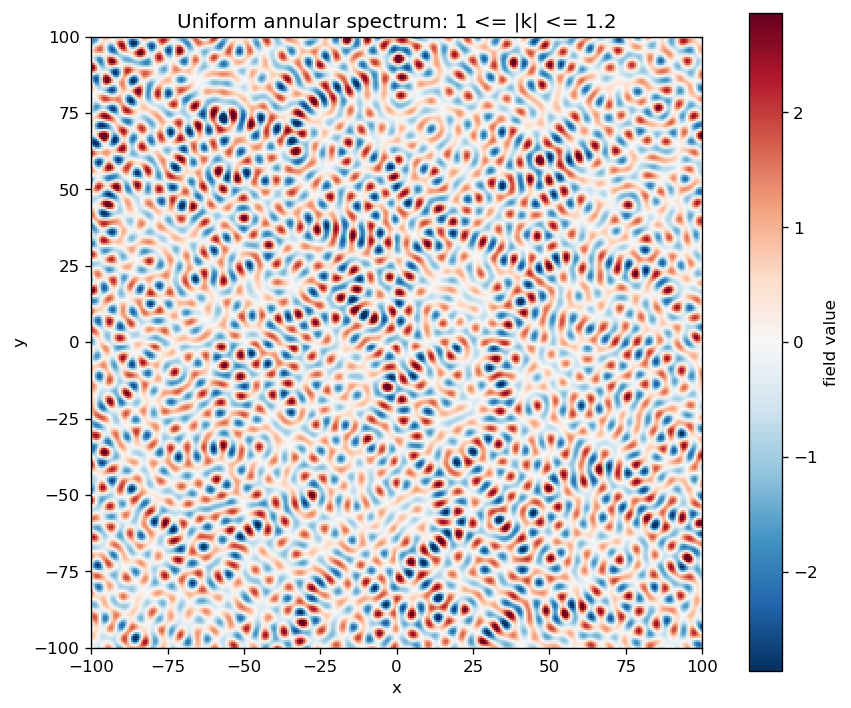}
    \includegraphics[width=0.49\linewidth]{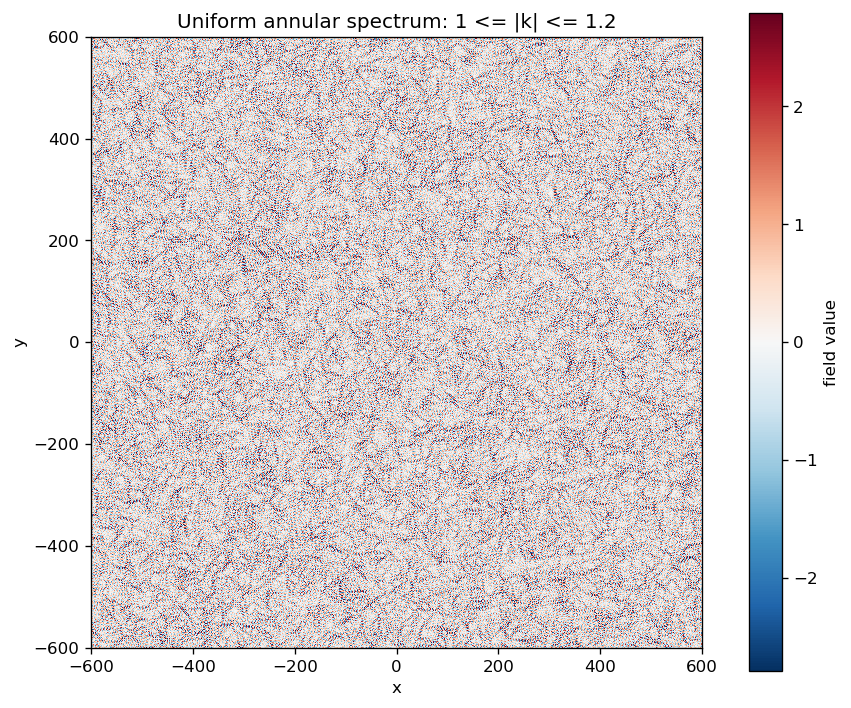}
    \caption{Gaussian field sample corresponding to the spectral measure which is uniform in an annulus $1\le |\zeta|\le 1.2$, in $200 \times 200$ and $1000 \times 1000$ boxes.}
    \label{fig:annular-fields}
\end{figure}

Finally, we want to make a comment about the filaments. 
It is known that they are a very special property of the RPW. It seems to be related to the fact that it is a monochromatic wave.  In particular, if we allow a range of different frequencies, then the filaments disappear. This is consistent with the physics paper observation in \cite{oconnor_gehlen_heller_1987}  that the string-like structures visible in monochromatic random waves are destroyed when the field is formed from a range of wavelengths. Hence, slow or oscillatory covariance decay by itself is not the decisive criterion for the existence of filaments. 

\pgap 

Having said that, we want to note that the exact monochromaticity is not absolutely necessary, which is pointed out in the next remark. 

\begin{remark}
Consider the field which is the sum of two RPW with different frequencies. Its spectral measure is the sum of arclengths on two circles. The proof of Theorem \ref{thm:main-thm} will go through almost verbatim in this case. The filaments of two RPW do not cancel out and can be seen in the sum. 
\end{remark} 

\pgap 

In Figure \ref{fig:annular-fields}, we can see a sample of a limited band wave where the spectral measure is the uniform measure in a thin annulus. Since this spectral measure is close to the uniform measure on the unit circle, on small scales, this field looks like the RPW. In particular, on small scales (left) one can see filaments. On the larger scales, these fields look different, and one can see short filaments that do not extend beyond a certain scale. We also recommend watching the `Worms movie' \cite{Warms_movie}. It is an animation which shows how a field changes when the annulus gets thinner and converges to the unit circle. 

This illustrates the distinction between local appearance and large-scale fluctuation behaviour: a field may display short filament-like patterns on moderate scales because its spectrum is concentrated near a circle, while still belonging to a different large-scale universality class.

To illustrate it with a short calculation, consider the isotropic band-limited field whose spectral density is uniform on the annulus
\[
    A_{a,b}:=\{\xi\in\bR^2:a\leq |\xi|\leq b\}, \qquad 0\leq a<b.
\]
Writing
\[
    \sigma_{a,b}(\xi)
    =\frac{1}{\pi(b^2-a^2)}\bI_{A_{a,b}}(\xi),
\]
the covariance is
\begin{align*}
    K_{a,b}(r)
    &=\int_{\bR^2}e^{i x\cdot \xi}\sigma_{a,b}(\xi)\td \xi \\
    &=\frac{2}{b^2-a^2}\int_a^b J_0(sr)s\td s
      =\frac{2}{(b^2-a^2)r}\big(bJ_1(br)-aJ_1(ar)\big),
      \qquad r=|x|.
\end{align*}
The Bessel asymptotic gives $K_{a,b}(r)=O(r^{-3/2})$ as $r\to\infty$, hence $K_{a,b} \notin L^1(\bR^2)$ i.e. long-range field just as RPW.
More importantly, $\sigma_{a,b}$ is an $L^\infty$ density, and Young's
inequality gives
\[
    \|\sigma_{a,b}^{*q}\|_\infty
    \leq \|\sigma_{a,b}\|_\infty \|\sigma_{a,b}\|_1^{q-1}
    =\|\sigma_{a,b}\|_\infty
    \qquad (q\geq1).
\]
Thus, all convolution powers\footnote{The importance of the convolution powers will become apparent later} have bounded density at the origin. The Riesz
singularity is absent, so the leading chaoses of smooth local functionals should have $R^2$-order local covariance
and white-noise scaling, rather than the $R^3$ fractional-Gaussian scaling of
Theorem \ref{thm:main-thm} in dimension $2$.

\subsection{Fields with a spectral singularity at the origin}

We now consider a class of smooth stationary Gaussian fields that have a spectral measure in a neighbourhood of the origin but the spectral density blows up. We will see that in this case, we still have long-range correlations but they belong to different universality classes. For simplicity and to illustrate the point, we are going to consider the simplest case. 

Fix
\[
    0<\alpha<2
\]
and let $f_\alpha:\bR^2\to\bR$ be the centered stationary Gaussian field with
spectral density
\begin{equation} \label{eq:smooth-gf-spectral-density}
    \sigma_\alpha(\xi)
    =
    c_\alpha |\xi|^{-\alpha}\bI_{\{|\xi|<1\}},
    \qquad
    c_\alpha=\frac{2-\alpha}{2\pi}.
\end{equation}
Thus $\int_{\bR^2}\sigma_\alpha(\xi)\td \xi=1$, so $f_\alpha(x)$ has unit
variance, and
\[
    K_\alpha(x-y)
    :=\bE[f_\alpha(x)f_\alpha(y)]
    =
    \int_{\bR^2}e^{i(x-y)\cdot\xi}\sigma_\alpha(\xi)\td \xi .
\]
The condition $\alpha<2$ is exactly the local integrability of
$|\xi|^{-\alpha}$ at the origin in dimension two. Since the spectral measure
has compact support, all spectral moments are finite; consequently
$f_\alpha$ has a smooth version.

\begin{figure}[ht]
    \centering
    \includegraphics[width=0.9\linewidth]{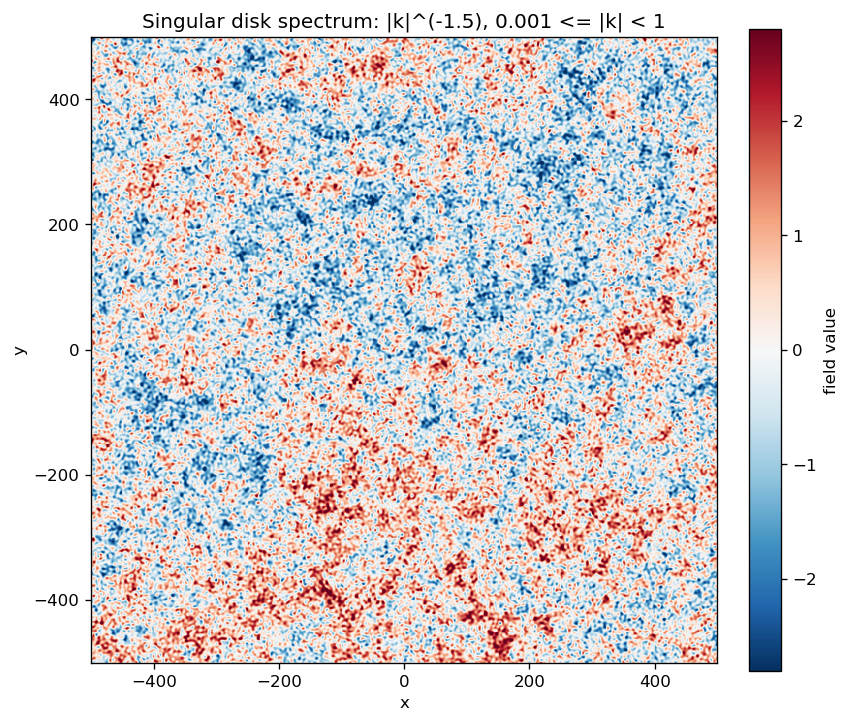}
    \caption{A sample of the smooth field $f_\alpha$, with $\alpha=1.5$}
    \label{fig:SGF-alpha}
\end{figure}

\begin{theorem} \label{thm:smooth-low-frequency-gf}
Let $0<\alpha<2$ and let $f_\alpha$ be the Gaussian field with spectral
density \eqref{eq:smooth-gf-spectral-density}. Consider a linear functional of the excursion set $A^{(\alpha)}_{R,u}$ defined in the same way as before by \eqref{eq:area-statistics} (with $f_\alpha$ instead of $F_d$). Then for every fixed
$u\in\bR$
\[
    R^{-1-\alpha/2}
    \left(A^{(\alpha)}_{R,u}
    -\bE[A^{(\alpha)}_{R,u}]\right)
    \Rightarrow
    \gamma(u)\sqrt{2\pi(2-\alpha)}\,\FGF_{\alpha/2}(\bR^2), \quad R\to \infty
\]
in distribution in $\mathcal S'(\bR^2)$, where, as before, $\gamma(u)$ is the Gaussian density.
\end{theorem}

\begin{remark}
The theorem is stated in the genuinely singular low-frequency regime
$0<\alpha<2$. At $\alpha \leq 0$ the spectral density is bounded at the origin and
the first chaos no longer separates from the higher chaoses by scale; the
natural limit is then of white-noise type rather than the positive-parameter
FGF above. Compare this with \cite[Theorem 2.4.9]{gass2026}.
\end{remark}

\begin{remark}
Note that for $ \alpha=1$ the formula is the same up to the constant factor as in the RPW case but the field $f_1$ is very different from the RPW.    
\end{remark}

\subsection{Strategy of the proof}
\label{ss: strategy}

Both Theorems \ref{thm:main-thm}  and \ref{thm:smooth-low-frequency-gf} follow the same general scheme. One expands the centred excursion indicator into Hermite chaoses, identifies the lowest chaos whose spectral measure has the relevant singular behaviour at the origin, proves that this chaos gives the leading contribution at the fluctuation scale, and then applies a spectral central limit theorem to obtain a fractional Gaussian limit in $\mathcal{S}'(\bR^d)$. The difference between the two settings lies in which chaos dominates.  For Berry’s random wave, there is no spectral density near the origin; hence the first chaos is negligible, and the second chaos determines the limit.  For the planar fields with spectral blow-up, the singularity is already visible in the first chaos.

The common strategy is the following: 
Let $\{H_q\}_{q\geq0}$ be the probabilists' Hermite polynomials and set
\[
    g_u(z):=\bI[z>u]-\bP(Z>u),
    \qquad Z\sim N(0,1).
\]
Then
\begin{equation} \label{eq:monochromatic-indicator-chaos}
    g_u(z)=\sum_{q\geq1}a_q(u)H_q(z),
    \qquad
    a_q(u)=\gamma(u)\frac{H_{q-1}(u)}{q!},
\end{equation}
with convergence in the standard Gaussian $L^2$ space. Hence
\[
    a_1(u)=\gamma(u),
    \qquad
    a_2(u)=\frac{u\gamma(u)}{2}.
\]
For $q\geq1$ define
\[
    Z_{q,R}(\varphi)
    :=
    \int_{\bR^d}H_q(f(x))\varphi(x/R)\td x,
\]
where $f$ is the field that we are studying.  Then
\begin{equation} 
    A_{R,u}(\varphi)-\bE[A_{R,u}(\varphi)]
    =
    \sum_{q\geq1}a_q(u)Z_{q,R}(\varphi)
\end{equation}
in $L^2(\Omega)$. Here $A$ is the relevant linear functional as in Theorems \ref{thm:main-thm} and \ref{thm:smooth-low-frequency-gf}, and $\phi$ is a test function.

A simple computation will show that 
\begin{equation} 
    \cov\big(Z_{q,R}(\varphi),Z_{q,R}(\psi)\big)
    =
    q!R^d\int_{\bR^d}
    p_{q}(\xi/R)\widehat\varphi(\xi)\widehat\psi(-\xi)\td \xi ,
\end{equation}
where $p_q$ is the density of the convolution of $q$ copies of the spectral measure of $f$ (if there is no density, this should be interpreted as the integral against the corresponding measure). 

\pgap 

This is exactly where the proofs of two theorems start to diverge, since these convolutions could have very different behaviour depending on the spectral measure. After that, we identify the leading terms in the chaos expansion, argue that the tail does not contribute, and finally use a special version of the CLT to show convergence to FGF.

\subsection{Possible connection to GFF and SLE}
Nothing in this section is rigorous. One should treat it as a heuristic or a far-fetched, vague conjecture. The most famous example of a fractional Gaussian field is the two-dimensional Gaussian free field (GFF). In
our notation
\[
    h_{\mathrm{GFF}}=\mathrm{FGF}_1(\bR^2),
    \qquad H=1-\frac{2}{2}=0 .
\]

This field is also a random distribution and has no level sets, but there is a way to make sense out of the nodal lines of GFF and show that they are (ignoring a lot of technicalities) given by Schramm-Loewner Evolution SLE($4$) curves. In particular, a global nodal line of a discrete GFF converges to SLE($4$) \cite{schramm_sheffield_contour_2009}.

\pgap 

For the field appearing in our scaling limit,
\[
    X=\mathrm{FGF}_{1/2}(\bR^2),\qquad
    \mathbf E[X(\varphi)X(\psi)]
    =
    \frac1{(2\pi)^2}\int_{\bR^2}
    \widehat\varphi(\xi)\widehat\psi(-\xi)|\xi|^{-1}\td \xi ,
\]
the Hurst exponent is
\[
    H=s-\frac d2=\frac{1}{2}-1=-\frac{1}{2}.
\]

The results about GFF and SLE($4$) rely heavily on the fact that GFF has a domain Markov property.
Lodhia, Sheffield, Sun and Watson formulate the corresponding domain Markov
analogue for fractional Gaussian fields by orthogonal projection in the
$\dot H^s$ Cameron--Martin space \cite[Sections~4--5]{lodhia_fractional_2016}.
If $D\subset\bR^2$ is an allowable domain, then conditionally on the
exterior field the restriction to $D$ decomposes as
\[
    X|_D=X_D^{\mathrm{har}}+X_D^0,
\]
where $X_D^0$ is an independent zero-boundary $\mathrm{FGF}_{1/2}$ on $D$
and $X_D^{\mathrm{har}}$ is the $1/2$-harmonic extension of the exterior
data, i.e. $(-\Delta)^{1/2}X_D^{\mathrm{har}}=0$ in $D$ in the distributional
sense. This is not the local Markov property of the GFF: since
$(-\Delta)^{1/2}$ is nonlocal, the conditional mean depends on the full
exterior data, not only on infinitesimal boundary values. The local-set
formalism nevertheless extends to the FGF in this projection sense
\cite[Section~8]{lodhia_fractional_2016}. This suggests that the analogy might go even further.

\pgap

The physics prediction for the interface exponent is consistent with this
identification of $H$. Strictly speaking, $\mathrm{FGF}_{1/2}(\bR^2)$ is a
random distribution, so its pointwise level sets are not defined without a
regularization. De Castro, Lukovi\'c, Pompanin, Andrade and Herrmann study in \cite{decastro_schramm_2018}
the corresponding lattice/regularized correlated Gaussian surfaces whose spectral density satisfies
\[
    m(\xi)\sim |\xi|^{-\beta},\qquad \beta=2(H+1), 
\]
Since a two-dimensional $\mathrm{FGF}_s$ has
spectral density proportional to $|\xi|^{-2s}$ and $H=s-1$, this is the same
spectral scaling:
\[
    2s=2(H+1).
\]
For $s=1/2$ one has $\beta=1$ and $H=-1/2$. The numerical SLE prediction in \cite{decastro_schramm_2018} uses the
complete perimeter fractal-dimension law
\[
    d_f(H)=\frac32-\frac H3,
    \qquad -\frac34\leq H\leq 0,
\]
conjectured in \cite{schrenk_percolation_2013}, together with the SLE
dimension formula
\[
    d_f=1+\frac{\kappa}{8},\qquad \kappa\leq 8.
\]
Thus
\[
    \kappa(H)=8(d_f(H)-1)
    =4-\frac{8H}{3}.
\]
Substituting $H=-1/2$ gives
\[
    \kappa=4+\frac43=\frac{16}{3}.
\]
In the same simulations, the Loewner driving function has variance
$\langle \zeta_t^2\rangle\sim \kappa t$, negligible time correlations after
the lattice-scale transient, and approximately Gaussian one-time marginals.
Thus the expected scaling limit of the regularized iso-height interface with
the $\mathrm{FGF}_{1/2}(\bR^2)$ spectral exponent is $\mathrm{SLE}_{16/3}$,
in the numerical/universality sense of \cite{decastro_schramm_2018}. This
should be distinguished from the rigorous Schramm--Sheffield GFF theorem above.
The same value $\kappa=16/3$ is the rigorous scaling-limit parameter for
critical FK-Ising interfaces, while spin-Ising interfaces have parameter
$3$ \cite{chelkak_ising_interfaces_2014}.

\pgap 

It is also interesting to compare the nodal lines of smooth fields like RPW or even better Bargmann-Fock field\footnote{This is a field with covariance given by the Gaussian kernel.}. The Bogomolny-Schmit conjecture states that their nodal lines at large scales look like the critical percolation interfaces and thus converge to SLE($6$) curves. Importantly, these two statements could be true (at least conjecturally) simultaneously. For example, consider a covariance kernel which is fast-decaying at infinity and has logarithmic divergence at the origin. Let us mollify the field; then it will be a smooth, symmetric field with fast-decaying correlations. The generalised Bogomolny-Schmit conjecture suggests that the scaling limit is in the universality class of SLE($6$). On the other hand, as the degree of mollification goes to zero, we recover the logarithmic singularity, and hence the nodal lines should converge to SLE($4$). There is no contradiction since we consider different limits. 

\section{Proof of Theorem \ref{thm:main-thm}}
We now implement the general strategy described in Section \ref{ss: strategy} in the case of the Berry's random wave $F_d$.
\subsection*{Preliminary computation: the two-step spherical density.}
In this section, we compute the density of the convolution square of the spectral measure of the RPW. We provide this for the sake of completeness since this is not a new result. See, for example, \cite[Lemma 2.2]{grotto_fluctuations_2024}. 

\pgap 

Let $U,V$ be independent with law $\mu_d$. The density of
$T:=U\cdot V$ on $[-1,1]$ is
\[
    b_d(1-t^2)^{(d-3)/2},
    \qquad
    b_d:=\frac{\Gamma(d/2)}
    {\sqrt{\pi}\Gamma((d-1)/2)}.
\]
Indeed, by rotation invariance we may condition on $U=e_d$. Then
$T$ has the same law as the last coordinate $V_d$ of a uniform point on
$S^{d-1}$. The slice $\{v\in S^{d-1}:v_d=t\}$ is a copy of
$S^{d-2}$ with radius $(1-t^2)^{1/2}$, so the coarea formula gives
\[
    \bP(V_d\in\td t)
    =
    \frac{\omega_{d-2}}{\omega_{d-1}}
    (1-t^2)^{(d-3)/2}\td t
    =
    b_d(1-t^2)^{(d-3)/2}\td t.
\]
Since $|U+V|=\sqrt{2+2T}$, the radial density of $|U+V|$ is
\[
    b_d r^{d-2}\left(1-\frac{r^2}{4}\right)^{(d-3)/2},
    \qquad 0<r<2.
\]
Dividing by the Euclidean surface factor $\omega_{d-1}r^{d-1}$ gives the density 
\begin{equation} \label{eq:monochromatic-p2-density}
    p_{d,2}(\xi)
    =
    \kappa_d\frac{1}{|\xi|}
    \left(1-\frac{|\xi|^2}{4}\right)^{(d-3)/2}
    \bI_{\{0<|\xi|<2\}},
    \qquad
    \kappa_d:=
    \frac{\Gamma(d/2)^2}
    {2\pi^{(d+1)/2}\Gamma((d-1)/2)}.
\end{equation}
In particular $p_{d,2}(\xi)\sim \kappa_d|\xi|^{-1}$ at the origin. For
$d=2$ this is
\[
    p_{2,2}(\xi)=
    \frac{1}{\pi^2|\xi|\sqrt{4-|\xi|^2}}\bI_{\{0<|\xi|<2\}},
\]
and for $d=3$ it is
\[
    p_{3,2}(\xi)=
    \frac{1}{8\pi|\xi|}\bI_{\{0<|\xi|<2\}}.
\]

\subsection*{Wiener chaos expansion and covariance asymptotics.}

As explained in Section \ref{ss: strategy}, we start by expanding the indicator function in Hermite polynomials \eqref{eq:monochromatic-indicator-chaos}.

For $q\geq1$ define
\[
    Z^{(d)}_{q,R}(\varphi)
    :=
    \int_{\bR^d}H_q(F_d(x))\varphi(x/R)\td x.
\]
Then
\begin{equation} \label{eq:monochromatic-area-chaos}
    A^{(d)}_{R,u}(\varphi)-\bE[A^{(d)}_{R,u}(\varphi)]
    =
    \sum_{q\geq1}a_q(u)Z^{(d)}_{q,R}(\varphi)
\end{equation}
in $L^2(\Omega)$. If $X,Y$ are standard Gaussian variables with correlation $\rho$, then the
generating function
\[
    e^{tz-t^2/2}=\sum_{n\geq 0}H_n(z)\frac{t^n}{n!}
\]
gives
\[
    \bE\left[e^{tX-t^2/2}e^{sY-s^2/2}\right]
    =e^{\rho ts}
    =\sum_{n\geq 0}\rho^n\frac{(ts)^n}{n!}.
\]
Comparing the coefficients of $t^q s^q$ yields
\[
    \bE[H_q(X)H_q(Y)]=q!\rho^q.
\] 

Therefore
\begin{equation} \label{eq:monochromatic-chaos-covariance}
    \cov\big(Z^{(d)}_{q,R}(\varphi),Z^{(d)}_{q,R}(\psi)\big)
    =
    q!\iint_{\bR^d\times\bR^d}
    K_d(x-y)^q\varphi(x/R)\psi(y/R)\td x\td y.
\end{equation}
Since
\[
    K_d(z)^q
    =
    \int_{\bR^d}e^{iz\cdot\xi}\td\mu_d^{*q}(\xi),
\]
whenever $\mu_d^{*q}$ has Lebesgue density $p_{d,q}$,
\begin{equation} \label{eq:monochromatic-random-walk-covariance}
    \cov\big(Z^{(d)}_{q,R}(\varphi),Z^{(d)}_{q,R}(\psi)\big)
    =
    q!R^d\int_{\bR^d}
    p_{d,q}(\xi/R)\widehat\varphi(\xi)\widehat\psi(-\xi)\td \xi .
\end{equation}

For $q=1$, the spectral measure is supported on $S^{d-1}$, away from the
origin. Thus for every $N$,
\begin{equation} \label{eq:monochromatic-first-chaos}
    \cov\big(Z^{(d)}_{1,R}(\varphi),Z^{(d)}_{1,R}(\psi)\big)
    =
    O_{N,\varphi,\psi}(R^{-N}).
\end{equation}
For $q=2$, \eqref{eq:monochromatic-p2-density} and
\eqref{eq:monochromatic-random-walk-covariance} give
\begin{align}
    &\cov\big(Z^{(d)}_{2,R}(\varphi),Z^{(d)}_{2,R}(\psi)\big)
    \notag \\
    &\qquad
    =
    2\kappa_d R^{d+1}
    \int_{|\xi|<2R}
    \frac{\widehat\varphi(\xi)\widehat\psi(-\xi)}{|\xi|}
    \left(1-\frac{|\xi|^2}{4R^2}\right)^{(d-3)/2}
    \td \xi. \label{eq:monochromatic-second-chaos-exact}
\end{align}
The factor in parentheses converges to $1$ locally uniformly. The part
$|\xi|>R$ is $o(1)$ after division by $R^{d+1}$, by the rapid decay of
$\widehat\varphi$ and $\widehat\psi$; when $d=2$ the possible singularity at
$|\xi|=2R$ is integrable and gives the same conclusion in polar coordinates.
Hence
\begin{equation} \label{eq:monochromatic-second-chaos-covariance}
    R^{-(d+1)}
    \cov\big(Z^{(d)}_{2,R}(\varphi),Z^{(d)}_{2,R}(\psi)\big)
    \to
    2\kappa_d
    \int_{\bR^d}
    \frac{\widehat\varphi(\xi)\widehat\psi(-\xi)}{|\xi|}\td \xi .
\end{equation}

We next bound the remaining chaoses. The covariance kernel has the Bessel
form
\[
    K_d(z)=
    2^{d/2-1}\Gamma(d/2)
    \frac{J_{d/2-1}(|z|)}{|z|^{d/2-1}},
\]
with the value at $z=0$ understood by continuity. Consequently, by the
standard Bessel asymptotic estimate, in the form used for Euclidean random
waves in \cite[Lemma~2.5]{grotto_fluctuations_2024},
\begin{equation} \label{eq:monochromatic-kernel-decay}
    |K_d(z)|\leq C_d(1+|z|)^{-(d-1)/2},
    \qquad z\in\bR^d.
\end{equation}
For Schwartz functions,
\[
    \int_{\bR^d}
    |\varphi((y+z)/R)\psi(y/R)|\td y
    \leq
    C_{\varphi,\psi,M}R^d(1+|z|/R)^{-M}
\]
for every $M$. Combining this with \eqref{eq:monochromatic-kernel-decay},
for fixed $q\geq3$ we obtain
\begin{equation} \label{eq:monochromatic-fixed-higher-chaos-bound}
    \left|
    \cov\big(Z^{(d)}_{q,R}(\varphi),Z^{(d)}_{q,R}(\psi)\big)
    \right|
    \leq
    C_{q,\varphi,\psi}R^d B_{d,q}(R),
\end{equation}
where, writing $\alpha_{d,q}=q(d-1)/2$,
\[
    B_{d,q}(R)=
    \begin{cases}
        R^{d-\alpha_{d,q}},& \alpha_{d,q}<d,\\
        \log R,& \alpha_{d,q}=d,\\
        1,& \alpha_{d,q}>d.
    \end{cases}
\]
To see this, write $z=x-y$ in
\eqref{eq:monochromatic-chaos-covariance}. The preceding two estimates give
\[
\begin{aligned}
    \left|
    \cov\big(Z^{(d)}_{q,R}(\varphi),Z^{(d)}_{q,R}(\psi)\big)
    \right|
    &\leq
    C_{q,\varphi,\psi,M}R^d
    \int_{\bR^d}
    (1+|z|)^{-\alpha_{d,q}}(1+|z|/R)^{-M}\td z .
\end{aligned}
\]
Taking $M>d+1$ and splitting the radial integral into $|z|\leq R$ and
$|z|>R$ gives exactly the three alternatives defining $B_{d,q}(R)$. This is
the smooth-weight analogue of the Euclidean polyspectra estimate in
\cite[Theorem~1.3]{grotto_fluctuations_2024}; the Bessel input is the decay
estimate \cite[Lemma~2.5]{grotto_fluctuations_2024}, which is written there
for the same Berry covariance kernel.

Since $q\geq3$ and $d\geq2$ imply $\alpha_{d,q}>d-1$, we have
$B_{d,q}(R)=o(R)$. Thus every fixed chaos $q\geq3$ is
$o(R^{d+1})$ at the variance scale.

\pgap 

It remains to control the infinite tail uniformly in $q$. Choose an integer
$Q=Q(d)\geq3$ such that $Q(d-1)/2>d$. Then $|K_d|^Q\in L^1(\bR^d)$, and
$|K_d|\leq1$ because $K_d$ is a covariance kernel. Hence for every $q\geq Q$,
\[
    \left|
    \iint
    K_d(x-y)^q\varphi(x/R)\psi(y/R)\td x\td y
    \right|
    \leq
    C_{\varphi,\psi}R^d .
\]
Since $\sum_{q\geq1}q!a_q(u)^2=\bE[g_u(Z)^2]<\infty$, the contribution of
all $q\geq Q$ to the variance of
\eqref{eq:monochromatic-area-chaos} is $O_{\varphi}(R^d)$.
Together with \eqref{eq:monochromatic-first-chaos} and
\eqref{eq:monochromatic-fixed-higher-chaos-bound}, this proves
\begin{equation} \label{eq:monochromatic-remainder-negligible}
    R^{-(d+1)}
    \bE\left[
        \left|
        \sum_{q\neq2}a_q(u)Z^{(d)}_{q,R}(\varphi)
        \right|^2
    \right]
    \to0.
\end{equation}

\medskip

\subsection*{Gaussian limit of the second chaos.}
We use the following rank-two smooth-weight form of the Maini--Nourdin
spectral CLT \cite[Theorem~2]{maini_spectral_2024}. Let $B$ be a continuous
stationary isotropic unit-variance Gaussian field, and let $m_2(\xi)\td\xi$
be the spectral measure of $\bE[H_2(B(0))H_2(B(\cdot))]$. For
\[
    Y_L(\phi):=\int_{\bR^d}H_2(B(x))\phi(x/L)\td x,
    \qquad \phi\in\mathcal S(\bR^d),
\]
assume that, for some locally finite non-atomic measure $\nu$,
\begin{equation} \label{eq:mn-rank-two-input}
    L^{-1}m_2(\xi/L)\td \xi
    \xrightarrow[L\to\infty]{\mathcal S'}
    \nu(\td \xi).
\end{equation}
Then, for every $\phi_1,\ldots,\phi_m\in\mathcal S(\bR^d)$,
\begin{equation} \label{eq:mn-rank-two-output}
    \left(L^{-(d+1)/2}Y_L(\phi_j)\right)_{j=1}^m
    \Rightarrow
    \left(G(\phi_j)\right)_{j=1}^m,
    \qquad
    \bE[G(\phi)G(\psi)]
    =
    \int_{\bR^d}\widehat\phi(\xi)\widehat\psi(-\xi)\nu(\td\xi).
\end{equation}
For a Schwartz weight this is obtained from the compact-domain statement of
\cite[Theorem~2]{maini_spectral_2024} by smooth truncation; the covariance
identity above controls the truncation error.

We now verify \eqref{eq:mn-rank-two-input}. In the present field,
\[
    \bE[H_2(F_d(0))H_2(F_d(z))]=2K_d(z)^2,
    \qquad
    m_2(\xi)=2p_{d,2}(\xi).
\]
Hence \eqref{eq:monochromatic-p2-density} gives
\begin{equation} \label{eq:monochromatic-rank-two-spectral-limit}
    R^{-1}m_2(\xi/R)\td \xi
    =
    2\kappa_d
    \frac{\bI_{\{|\xi|<2R\}}}{|\xi|}
    \left(1-\frac{|\xi|^2}{4R^2}\right)^{(d-3)/2}
    \td \xi
    \xrightarrow[R\to\infty]{\mathcal S'}
    2\kappa_d\frac{\td \xi}{|\xi|}.
\end{equation}
The convergence is exactly the limit in
\eqref{eq:monochromatic-second-chaos-covariance}, tested against
$\widehat\phi(\xi)\widehat\psi(-\xi)$. The limit measure is locally finite
and non-atomic because $|\xi|^{-1}\in L^1_{\operatorname{loc}}(\bR^d)$ for
$d\geq2$. Applying \eqref{eq:mn-rank-two-output} gives, for every
$\varphi_1,\ldots,\varphi_m\in\mathcal S(\bR^d)$,
\begin{equation} \label{eq:monochromatic-second-chaos-clt}
    \left(
        R^{-(d+1)/2}Z^{(d)}_{2,R}(\varphi_j)
    \right)_{j=1}^m
    \Rightarrow
    \left(G_d(\varphi_j)\right)_{j=1}^m,
\end{equation}
where $G_d$ is centered Gaussian with covariance
\begin{equation} \label{eq:monochromatic-gd-covariance}
    \bE[G_d(\varphi)G_d(\psi)]
    =
    2\kappa_d
    \int_{\bR^d}
    \frac{\widehat\varphi(\xi)\widehat\psi(-\xi)}{|\xi|}\td \xi .
\end{equation}

\medskip

\subsection*{Identification of the limit.}
By \eqref{eq:monochromatic-remainder-negligible},
\[
    R^{-(d+1)/2}
    \big(A^{(d)}_{R,u}(\varphi)-\bE[A^{(d)}_{R,u}(\varphi)]\big)
    -
    a_2(u)R^{-(d+1)/2}Z^{(d)}_{2,R}(\varphi)
    \to0
\]
in $L^2(\Omega)$. Combining this with
\eqref{eq:monochromatic-second-chaos-clt}, the finite-dimensional limit of
the centered excursion field is Gaussian with covariance
\[
    a_2(u)^2\,2\kappa_d
    \int_{\bR^d}
    \frac{\widehat\varphi(\xi)\widehat\psi(-\xi)}{|\xi|}\td \xi
    =
    \frac{u^2\gamma(u)^2\kappa_d}{2}
    \int_{\bR^d}
    \frac{\widehat\varphi(\xi)\widehat\psi(-\xi)}{|\xi|}\td \xi .
\]
Comparing this with \eqref{eq:monochromatic-fgf-covariance}, the multiplier
in front of $h_d$ must satisfy
\[
    c_{d,u}^2
    =
    (2\pi)^d\frac{u^2\gamma(u)^2\kappa_d}{2}.
\]
Substituting the value of $\kappa_d$ from
\eqref{eq:monochromatic-p2-density} gives exactly the stated constant
$c_{d,u}$.

\pgap 

Finally, the convergence holds in $\mathcal S'(\bR^d)$ by L{\'e}vy's
continuity theorem for generalized random fields. Indeed, for each $R$,
\[
    \left|
    R^{-(d+1)/2}
    \big(A^{(d)}_{R,u}(\varphi)-\bE[A^{(d)}_{R,u}(\varphi)]\big)
    \right|
    \leq
    R^{(d-1)/2}\|\varphi\|_{L^1(\bR^d)},
\]
so the random functional is tempered, since the $L^1$ norm is controlled by Schwartz seminorms. To see this, consider
\[
    \|\varphi\|_{L^1(\bR^d)}
    \leq
    \left(\int_{\bR^d}(1+|x|^2)^{-d}\td x\right)
    \sup_{x\in\bR^d}(1+|x|^2)^d|\varphi(x)|.
\]
The last factor is one of the defining seminorms of $\mathcal S(\bR^d)$
(equivalently, it is bounded by finitely many of the standard seminorms
$\sup_x |x^\alpha \partial^\beta\varphi(x)|$). Therefore, for each $R$ and
each outcome, $X_R$ is a continuous linear functional on $\mathcal S(\bR^d)$.
 The limiting characteristic functional
is
\[
    L(\varphi)
    =
    \exp\left\{
    -\frac{c_{d,u}^2}{2(2\pi)^d}
    \int_{\bR^d}\frac{|\widehat\varphi(\xi)|^2}{|\xi|}\td \xi
    \right\}.
\]
If $\varphi_n\to0$ in $\mathcal S(\bR^d)$, then
$\widehat\varphi_n\to0$ in $\mathcal S(\bR^d)$ and
\[
    \int_{\bR^d}\frac{|\widehat\varphi_n(\xi)|^2}{|\xi|}\td \xi
    \leq
    C\sup_{|\xi|\leq1}|\widehat\varphi_n(\xi)|^2
    +
    C\sup_{\xi\in\bR^d}(1+|\xi|)^{d+2}
    |\widehat\varphi_n(\xi)|^2,
\]
where the second term uses the integrability of
$(1+|\xi|)^{-(d+2)}|\xi|^{-1}$ at infinity. Hence $L$ is continuous at the
origin. 

The following version of Levy's Theorem \ref{thm: Levy} shows immediately that finite dimensional distribution convergence upgrades to convergence in strong topology, which in turn implies convergence in weak topology.

\begin{theorem}[Theorem~2.3 and
Corollary~2.4 of \cite{bierme_generalized_2018}]
\label{thm: Levy}
Let $(X_n)_{n\ge 1}$ be a sequence of generalized random fields in $\bR^d$. If the characteristic functions $\mathcal{L}_{X_n}$ converge pointwise to a functional $\mathcal{L} : \mathcal{S} \longrightarrow \bC$ which is continuous at $0$, then
there exists a generalized random field $X$ such that $\mathcal{L}_X = \mathcal{L}$ and $X_n$ converges in
distribution to $X$ with respect to the strong topology.
\end{theorem} 
  
This completes the proof of Theorem \ref{thm:main-thm}. 

\begin{remark}
For $d=2$ one has $\kappa_2=1/(2\pi^2)$ and $c_{2,u}=u\gamma(u)$, recovering
the $R^{3/2}$ planar scaling. For $d=3$ one has $\kappa_3=1/(8\pi)$ and
$c_{3,u}=\pi u\gamma(u)/\sqrt2$, recovering the three-dimensional
$R^2$ scaling.
\end{remark}

\begin{remark}
The ingredients used above are standard in the recent polyspectra literature.
The relation between Wiener-chaos variances for Euclidean random waves and
uniform random-walk densities is developed in
\cite{grotto_fluctuations_2024}; in particular, their Euclidean variance
asymptotics and Bessel estimates are the ball-domain analogues of the
smooth-weight bounds used here. The Gaussian limit for the second chaos is an
application of the spectral central limit theorem of
\cite{maini_spectral_2024}, where the rank-two polyspectrum is precisely the
two-step spectral convolution. In the present monochromatic model this
convolution is the spherical random-walk density $p_{d,2}$ computed explicitly
in \eqref{eq:monochromatic-p2-density}.
\end{remark}


\section{Proof of Theorem \ref{thm:smooth-low-frequency-gf}}

We now implement the general strategy described in Section \ref{ss: strategy} in the case of the RPW field to $f_\alpha$ with the spectral blow-up of order $\alpha$.
We follow the same Wiener chaos decomposition as in the Berry random wave case. We define for $q\geq1$ d
\[
    Z^{(\alpha)}_{q,R}(\varphi)
    :=
    \int_{\bR^2}H_q(f_\alpha(x))\varphi(x/R)\td x 
\]
and
\begin{equation} \label{eq:smooth-gf-area-chaos}
    A^{(\alpha)}_{R,u}(\varphi)
    -\bE[A^{(\alpha)}_{R,u}(\varphi)]
    =
    \sum_{q\geq1}a_q(u)Z^{(\alpha)}_{q,R}(\varphi)
\end{equation}
in $L^2(\Omega)$.

Let
\[
    p_{q,\alpha}:=\sigma_\alpha^{*q}.
\]
Since
\[
    K_\alpha(z)^q
    =
    \int_{\bR^2}e^{iz\cdot\xi}p_{q,\alpha}(\xi)\td \xi ,
\]
the same scaled Plancherel identity used above gives
\begin{equation} \label{eq:smooth-gf-q-covariance}
    I^{(\alpha)}_{q,R}(\varphi,\psi)
    :=
    \bE[
    Z^{(\alpha)}_{q,R}(\varphi)
    Z^{(\alpha)}_{q,R}(\psi)]
    =
    q!R^2\int_{\bR^2}
    p_{q,\alpha}(\xi/R)
    \widehat\varphi(\xi)\widehat\psi(-\xi)\td \xi .
\end{equation}

The first chaos already sees the singularity of the spectral density at the
origin. From \eqref{eq:smooth-gf-spectral-density},
\begin{align}
    I^{(\alpha)}_{1,R}(\varphi,\psi)
    &=
    R^2\int_{\bR^2}
    \sigma_\alpha(\xi/R)
    \widehat\varphi(\xi)\widehat\psi(-\xi)\td \xi \notag\\
    &=
    c_\alpha R^{2+\alpha}
    \int_{|\xi|<R}
    |\xi|^{-\alpha}
    \widehat\varphi(\xi)\widehat\psi(-\xi)\td \xi . \label{eq:smooth-gf-first-chaos}
\end{align}
Because $|\xi|^{-\alpha}$ is locally integrable and
$\widehat\varphi,\widehat\psi$ are rapidly decreasing,
\begin{equation} \label{eq:smooth-gf-first-chaos-limit}
    R^{-2-\alpha}I^{(\alpha)}_{1,R}(\varphi,\psi)
    \to
    c_\alpha
    \int_{\bR^2}
    |\xi|^{-\alpha}
    \widehat\varphi(\xi)\widehat\psi(-\xi)\td \xi .
\end{equation}

We now compare all higher chaoses with \eqref{eq:smooth-gf-first-chaos}. The
elementary Riesz-convolution estimate needed here is the following: for each
fixed $q\geq2$,
\begin{equation} \label{eq:smooth-gf-convolution-singularity}
    p_{q,\alpha}(\eta)
    =
    O\left(1+|\eta|^{-\beta_q}\right)
    \quad (\eta\to0),
    \qquad
    \beta_q:=\big(q\alpha-2(q-1)\big)_+,
\end{equation}
with the usual logarithmic replacement when $q\alpha=2(q-1)$. For example,
for $q=2$ this follows by splitting
\[
    \int |\zeta|^{-\alpha}|\eta-\zeta|^{-\alpha}
    \bI_{\{|\zeta|<1,\ |\eta-\zeta|<1\}}\td \zeta
\]
near $0$, near $\eta$, and away from both points. Iterating the same estimate
gives \eqref{eq:smooth-gf-convolution-singularity}. Since $\alpha<2$, for
every $q\geq2$ one has
\[
    \beta_q<\alpha .
\]
Therefore, for each fixed $q\geq2$,
\begin{equation} \label{eq:smooth-gf-fixed-higher-chaos}
    I^{(\alpha)}_{q,R}(\varphi,\psi)
    =
    O_{\varphi,\psi,q}(R^{2+\beta_q})
    +O_{\varphi,\psi,q}(R^2\log R)
    =
    o_{\varphi,\psi,q}(R^{2+\alpha}),
\end{equation}
where the logarithmic term is present only in the borderline case
$q\alpha=2(q-1)$.

It remains to control the infinite tail uniformly in $q$. The radial formula
\[
    K_\alpha(r)
    =
    (2-\alpha)\int_0^1 \rho^{1-\alpha}J_0(r\rho)\td \rho,
    \qquad r=|x|,
\]
and standard Bessel estimates give
\[
    |K_\alpha(r)|
    \leq C_\alpha(1+r)^{-\delta_\alpha},
    \qquad
    \delta_\alpha:=\min\{2-\alpha,3/2\}>0 .
\]
Choose $Q$ so large that $Q\delta_\alpha>2$. Then
$|K_\alpha|^Q\in L^1(\bR^2)$. Since $|K_\alpha|\leq1$, for all $q\geq Q$,
Young's inequality gives
\[
\begin{aligned}
    \left|
    \iint K_\alpha(x-y)^q
    \varphi(x/R)\psi(y/R)\td x\td y
    \right|
    &\leq
    \iint |K_\alpha(x-y)|^Q
    |\varphi(x/R)||\psi(y/R)|\td x\td y  \\
    &\leq
    \||K_\alpha|^Q\|_{L^1}
    \|\varphi(\cdot/R)\|_{L^2}
    \|\psi(\cdot/R)\|_{L^2} \\
    &=
    O_{\varphi,\psi}(R^2).
\end{aligned}
\]
Since $\sum_{q\geq1}q!a_q(u)^2=\bE[g_u(Z)^2]<\infty$, the tail
$q\geq Q$ contributes $O_{\varphi,u}(R^2)$ to the variance of
\eqref{eq:smooth-gf-area-chaos}. Combining this with
\eqref{eq:smooth-gf-fixed-higher-chaos}, we obtain
\begin{equation} \label{eq:smooth-gf-remainder-negligible}
    R^{-2-\alpha}
    \bE\left[
    \left|
    \sum_{q\geq2}a_q(u)Z^{(\alpha)}_{q,R}(\varphi)
    \right|^2
    \right]
    \to0 .
\end{equation}

Thus the centered excursion area is asymptotically equal, in $L^2$, to its
first chaotic projection:
\[
    R^{-1-\alpha/2}
    \left(A^{(\alpha)}_{R,u}(\varphi)
    -\bE[A^{(\alpha)}_{R,u}(\varphi)]\right)
    -
    \gamma(u)R^{-1-\alpha/2}Z^{(\alpha)}_{1,R}(\varphi)
    \to0
\]
in $L^2(\Omega)$. The remaining term is Gaussian for every $R$. By
\eqref{eq:smooth-gf-first-chaos-limit}, its limiting covariance is
\[
    \gamma(u)^2 c_\alpha
    \int_{\bR^2}
    |\xi|^{-\alpha}
    \widehat\varphi(\xi)\widehat\psi(-\xi)\td \xi .
\]
Comparing this with \eqref{eq:smooth-gf-spectral-density}, the limiting field is
\[
    \gamma(u)\,2\pi\sqrt{c_\alpha}\,h_{\alpha/2}
    =
    \gamma(u)\sqrt{2\pi(2-\alpha)}\,h_{\alpha/2}.
\]
This proves convergence of all finite-dimensional distributions, by
Cram{\'e}r--Wold and \eqref{eq:smooth-gf-remainder-negligible}.

Finally, the limiting characteristic functional is continuous at the origin
on $\mathcal S(\bR^2)$, because $|\xi|^{-\alpha}$ is locally integrable for
$\alpha<2$ and Schwartz decay controls infinity. L{\'e}vy's continuity
theorem for generalized random fields therefore upgrades the finite-dimensional
convergence to convergence in $\mathcal S'(\bR^2)$.

\appendix

\section{Fractional Gaussian Fields}
\label{ss: FGF}
This appendix about the notation and basic facts about fractional Gaussian
fields used in the paper. The reference is the survey of Lodhia, Sheffield,
Sun and Watson \cite{lodhia_fractional_2016}. We give definitions and
interpretations only, with no proofs.

\subsection{Whole-space definition}

Let $W$ be real white noise on $\bR^d$. Formally, the fractional Gaussian
field with parameter $s\in\bR$ is
\[
    \mathrm{FGF}_s(\bR^d)=(-\Delta)^{-s/2}W .
\]
Thus $\mathrm{FGF}_0$ is white noise, $\mathrm{FGF}_1$ is the Gaussian free
field, and $\mathrm{FGF}_2$ is the bi-Laplacian Gaussian field. The Hurst
parameter is
\[
    H=s-\frac d2 .
\]
This is the self-similarity parameter: if $h\sim\mathrm{FGF}_s(\bR^d)$,
then
\[
    h(a\cdot)\stackrel{\mathrm{law}}{=}a^Hh(\cdot),
    \qquad a>0,
\]
where the equality is interpreted distributionally. The rigorous meaning is as a centered Gaussian generalized field whose
covariance is the homogeneous Sobolev inner product
\begin{equation} \label{eq:fgf-fourier-covariance}
    \bE[h(\varphi)h(\psi)]
    =
    \langle\varphi,\psi\rangle_{\dot H^{-s}}
    =
    \frac{1}{(2\pi)^d}\int_{\bR^d}
    \widehat\varphi(\xi)\widehat\psi(-\xi)|\xi|^{-2s}\td\xi .
\end{equation}
Equivalently,
\[
    h(\varphi)=W((-\Delta)^{-s/2}\varphi),
\]
whenever the right-hand side is in $L^2$.

\subsection{Test functions and polynomials}

The singularity of $|\xi|^{-2s}$ at $\xi=0$ determines the correct test
space. If $s<d/2$, then $|\xi|^{-2s}$ is locally integrable and ordinary
Schwartz test functions are allowed:
\[
    \varphi,\psi\in\mathcal S(\bR^d).
\]
If $s\ge d/2$, then $H=s-d/2\ge0$ and the field is naturally defined modulo
polynomials of degree at most $\lfloor H\rfloor$. One tests against
\[
    \mathcal S_H(\bR^d)
    :=
    \left\{\varphi\in\mathcal S(\bR^d):
    \int_{\bR^d}\varphi(x)P(x)\td x=0
    \text{ for every polynomial }P,\ \deg P\le \lfloor H\rfloor
    \right\}.
\]
Equivalently,
\[
    \partial^\alpha\widehat\varphi(0)=0
    \qquad\text{for all }|\alpha|\le \lfloor H\rfloor .
\]
This is the same convention as Brownian motion modulo constants in one
dimension. For example, the two-dimensional GFF has $d=2$, $s=1$, and
$H=0$, so it is tested on mean-zero functions and is defined modulo an
additive constant.

\subsection{Covariance kernels}

The Fourier multiplier in \eqref{eq:fgf-fourier-covariance} is the cleanest
definition. In physical space one writes
\[
    \bE[h(\varphi)h(\psi)]
    =
    \int_{\bR^d}\int_{\bR^d}
    K_s(x,y)\varphi(x)\psi(y)\td x\td y ,
\]
where $K_s$ is a function or a distribution. If $0<s<d/2$, then
\[
    K_s(x,y)=C_{s,d}|x-y|^{2s-d}
    =
    C_{s,d}|x-y|^{2H},
\]
where $C_{s,d}>0$ is the Riesz-kernel constant in our Fourier convention. At the logarithmic value $s=d/2$, one obtains
\[
    K_{d/2}(x,y)=C_d\log\frac{1}{|x-y|}
\]
on mean-zero test functions. More generally, when $H$ is a nonnegative
integer, the kernel has the form
\[
    K_s(x,y)=C_{s,d}|x-y|^{2H}\log |x-y|
\]
again interpreted modulo the polynomial ambiguity described above.

\pgap 

For negative $s$, the FGF is a derivative of white noise. In particular,
if $s=-k$ with $k\in\bN\cup\{0\}$, then
\[
    K_{-k}(x,y)=(-\Delta)^k\delta(x-y),
\]
or equivalently
\[
    \bE[h(\varphi)h(\psi)]
    =
    \int_{\bR^d}\varphi(x)(-\Delta)^k\psi(x)\td x .
\]
For non-integer negative $s$, the covariance is the corresponding
fractional derivative of the delta distribution; the Fourier expression
\eqref{eq:fgf-fourier-covariance} is the safest notation.

\subsection{FGFs as long-range GFFs}

For $0<s<1$, the fractional Laplacian has the singular-integral form
\[
    (-\Delta)^s f(x)
    =
    c_{d,s}\,\mathrm{p.v.}\int_{\bR^d}
    \frac{f(x)-f(y)}{|x-y|^{d+2s}}\td y .
\]
This is nonlocal: points at all distances interact, with interaction kernel
of order $|x-y|^{-d-2s}$. Lodhia et al. construct a discrete
$\mathrm{FGF}_s$ in this range by replacing the nearest-neighbor gradient in
the discrete GFF energy with a fractional gradient. The resulting model is a
long-range Gaussian free field. Its potential theory is not Brownian
potential theory, as for the ordinary GFF, but the potential theory of an
isotropic $2s$-stable L{\'e}vy process.

\bibliographystyle{alpha}
\phantomsection
\addcontentsline{toc}{section}{References}
\bibliography{references}

\end{document}